\RequirePackage{fix-cm}
\documentclass[smallextended]{svjour3}       % onecolumn (second format)
\smartqed  % flush right qed marks, e.g. at end of proof
\usepackage{graphicx}
\usepackage{amssymb}
\usepackage{amsfonts}
\usepackage{amsmath}
\usepackage{epsfig}
\usepackage{psfrag}
\newcommand{\mR}{\mathbb{R}}

\newcommand{\bq}{\begin{equation}}
\newcommand{\eq}{\end{equation}}
\newcommand{\X}{\mathcal{X}}
\newcommand{\D}{\mathcal{D}}
\newcommand{\M}{\mathcal{M}}

\newcommand{\E}{\mathcal{E}}
\newcommand{\F}{\mathcal{F}}
\newcommand{\R}{\mathcal{R}}
\renewcommand{\L}{\mathcal{L}}
\begin{document}

\title{Dirac and Lagrange algebraic constraints in nonlinear port-Hamiltonian systems%\thanks{Grants or other notes
}

\titlerunning{Algebraic constraints in nonlinear port-Hamiltonian systems}        % if too long for running head

\author{Arjan van der Schaft \and \\Bernhard Maschke
}


\institute{A.J. van der Schaft \at
              Bernoulli institute for Mathematics, CS and AI, University of Groningen, the Netherlands \\
              Tel.: +31-50-3633731\\
              Fax: +31-50-3633800\\
              \email{a.j.van.der.schaft@rug.nl}           %  \\
          \and
           B. Maschke \at
           Universit\'e Claude Bernard Lyon-1, France
}

\date{Received: date / Accepted: date}

\maketitle

\begin{abstract}
After recalling standard nonlinear port-Hamiltonian systems and their algebraic constraint equations, called here Dirac algebraic constraints, an extended class of port-Hamiltonian systems is introduced. This is based on replacing the Hamiltonian function by a general Lagrangian submanifold of the cotangent bundle of the state space manifold, motivated by developments in \cite{barbero} and extending the linear theory as developed in \cite{DAE}, \cite{beattie}. The resulting new type of algebraic constraints equations are called Lagrange constraints. It is shown how Dirac algebraic constraints can be converted into Lagrange algebraic constraints by the introduction of extra state variables, and, conversely, how Lagrange algebraic constraints can be converted into Dirac algebraic constraints by the use of Morse families.

\end{abstract}

\section{Introduction}
\label{intro}
When modeling dynamical systems, the appearance of algebraic constraint equations next to differential equations is ubiquitous. This is especially clear in network modeling of physical systems, where the interconnections between different dynamical components of the overall system almost inevitably lead to algebraic constraints. On the other hand, the analysis and simulation of such systems of differential-algebraic equations (DAEs) poses delicate problems; especially in the nonlinear case, see e.g. the already classical treatise \cite{kunkel}, and the references quoted in there. The situation is even more challenging for control of DAE systems, and, up to now, most of the control literature is devoted to systems {\it without} algebraic constraints.

On the other hand, the DAE systems as resulting from the modeling of physical systems often have special properties, which makes their analysis, simulation and control potentially more feasible. A prominent framework for network modeling of multiphysics systems is {\it port-based modeling}, and the resulting theory of port-Hamiltonian systems; see e.g. \cite{maschkevdsbordeaux}, \cite{gsbm}, \cite{vanderschaftmaschkearchive}, \cite{dalsmo}, \cite{passivitybook}, \cite{geoplexbook}, \cite{NOW}. In \cite{phDAE}, see also \cite{passivitybook}, initial investigations were done on the algebraic constraint equations appearing in standard port-Hamiltonian systems; linear or nonlinear. The two main conclusions are that the algebraic constraints in such systems have index one once the Hamiltonian is non-degenerate in the state variables, and furthermore that generally elimination of algebraic constraints leads to systems that are again in port-Hamiltonian form. 

Very recently, see especially \cite{beattie}, {\it another} type of algebraic constraint equations in linear physical system models was studied. In \cite{DAE}, \cite{beattie} these were identified as arising from generalized port-based modeling with a state space that has higher dimension than the minimal number of energy variables, corresponding to {\it implicit} energy storage relations which can be formulated as Lagrangian subspaces. For linear time-varying systems this formulation has led to various interesting results on their index, regularization and numerical properties \cite{mmw}; see also \cite{bmd} for related developments.
The formulation of implicit energy storage relations is very similar to an independent line of research initiated in \cite{barbero}, where the Hamiltonian in nonlinear Hamiltonian dynamics is replaced by a general Lagrangian submanifold (motivated, among others, by optimal control).

The precise relation between the algebraic constraints as arising in linear standard port-Hamiltonian systems (called Dirac algebraic constraints) and those in linear generalized port-Hamiltonian systems with implicit storage (called Lagrange algebraic constraints) was studied in \cite{beattie}, \cite{DAE}. In particular, in \cite{DAE} it was shown how in this linear case Dirac algebraic constraints can be converted into Lagrange algebraic constraints, and conversely; in both cases by extending the state space (e.g., introduction of Lagrange multipliers).

In the present paper we will continue on \cite{DAE}, \cite{beattie}, by extending the theory and constructions to the nonlinear case; inspired by \cite{barbero}.

\section{The standard definition of port-Hamiltonian systems and Dirac algebraic constraints}
\label{sec:1}
\begin{figure}[t]
\begin{center}
\psfrag{D}[][]{\sf routing}
\psfrag{S}[][]{\sf storage}
\psfrag{R}[][]{\sf dissipation}
\psfrag{a}[][]{$e_S$}
\psfrag{b}[][]{$f_S$}
\psfrag{c}[][]{$e_R$}
\psfrag{d}[][]{$f_R$}
\psfrag{e}[][]{$e_P$}
\psfrag{f}[][]{$f_P$}
\includegraphics[width=12cm]{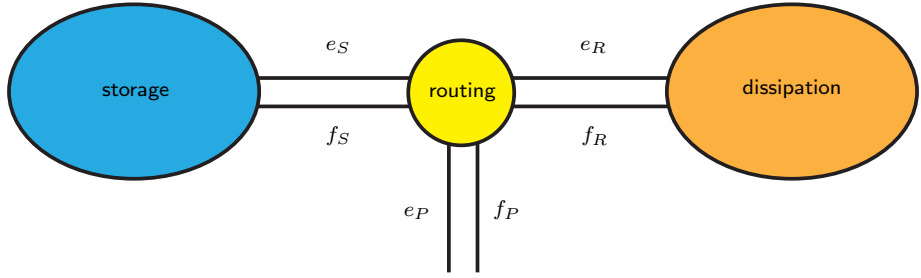}
\caption{From port-based modeling to port-Hamiltonian system.}
\label{fig:pHsystems}
\end{center}
\end{figure}

The standard definition of a port-Hamiltonian system, see e.g. \cite{vanderschaftmaschkearchive}, \cite{dalsmo}, \cite{phDAE}, \cite{passivitybook}, \cite{NOW} for more details and further ramifications, is depicted in Fig.~\ref{fig:pHsystems}. 

The system is modeled by identifying energy-storing elements $\mathcal{S}$ and energy-dissipating elements $\mathcal{R}$, which are linked to a central energy-routing structure, geometrically to be defined as a Dirac structure. This linking takes place via pairs $(f,e)$ of equally dimensioned vectors $f$ and $e$ (commonly called {\it flow} and {\it effort} variables). A pair $(f,e)$ of vectors of flow and effort variables defines a {\it port}, and the total set of variables $f,e$ is also called the set of {\it port variables}. 
Fig.~\ref{fig:pHsystems} shows three ports: the port $(f_S, e_S)$ linking to energy storage, the port $(f_R, e_R)$ corresponding to energy dissipation, and the external port $(f_P,e_P)$, by which the system interacts with its environment (including controller action). The scalar quantities $e^T_Sf_S$, $e^T_R f_R$, and   
 $e^T_P f_P$ denote the instantaneous powers transmitted through the links.
Any physical system that is represented (modeled) in this way defines a port-Hamiltonian system. Furthermore, experience has shown that even for very complex physical systems port-based modeling leads to satisfactory and insightful models; certainly for control purposes; see e.g. \cite{passivitybook}, \cite{gsbm}, \cite{NOW}, and the references quoted in there.

The definition of a {\it constant} Dirac structure starts from a finite-dimensional {\it linear space} of \textit{flows} $\mathcal{F}$ and the dual linear space of efforts $\mathcal{E}:= \mathcal{F}^*$. 
The \textit{power} $P$ on the total space $\mathcal{F}  \times \mathcal{E}$ of port variables is defined by the duality product
$P =\;< e \mid f >$. In the common case $\mathcal{F} \simeq \mathcal{E} \simeq \mR^k$ this simply amounts to $P= e^Tf$.
\begin{definition}\cite{Cour:90}, \cite{Dorf:93}
  \label{def:schaft_dn2.1}
  Consider a finite-dimensional linear space $\mathcal{F}$ with $\E = \F^*$. A subspace $\mathcal{D} \subset \mathcal{F} \times \mathcal{E}$ is a \textit{Dirac structure} if
  \begin{enumerate}
  \item[(i)] $< e \mid f > = 0$, for all $(f,e) \in \mathcal{D}$,
  \item[(ii)] $\dim \mathcal{D} = \dim \mathcal{F}$.
  \end{enumerate}
\end{definition}
The {\it maximal dimension} of any subspace $\mathcal{D} \subset \mathcal{F} \times \mathcal{E}$
satisfying the power-conservation property $(i)$ can be easily shown to be equal to $\dim \mathcal{F}$. Thus a constant Dirac structure is a {\it maximal} power-conserving subspace. 

In the nonlinear case, e.g. mechanical systems, we need to extend the notion of a constant Dirac structure to a {\it Dirac structure on a manifold}\footnote{All notions will be assumed to be smooth; i.e., $C^{\infty}$.} $\M$.
\begin{definition}\cite{Cour:90}, \cite{Dorf:93}
%\label{def:schaft_dn2.1}
A Dirac structure on a finite-dimensional manifold $\M$ is defined as a subbundle $\D \subset T\M \oplus T^*\M$ such that for every $x \in \M$ the subspace $\D(x) \subset T_x\M \times T^*_x\M$ is a constant Dirac structure as before.
\end{definition}
We note that in the standard definition \cite{Cour:90}, \cite{Dorf:93} of a Dirac structure on a manifold an additional {\it integrability} condition is imposed; generalizing the Jacobi identity for Poisson structures or closedness of symplectic forms. However, for many purposes this integrability condition is not needed, while on the other hand there are interesting port-Hamiltonian system classes (like mechanical systems with nonholonomic kinematic constraints) that do {\it not} satisfy this integrability condition \cite{passivitybook}, \cite{dalsmo}, \cite{NOW}.

The {\it dynamics} of port-Hamiltonian systems derives from the time-integration taking place in the energy storage. Let $f_S,e_S$ be the vector of flow and effort variables of the energy storage port. Time-integration of the flow vector $f_S$ leads to the equally dimensioned vector of state variables $x \in \mathcal{X}$ satisfying $\dot{x}=-f_S$. Energy storage in a standard port-Hamiltonian system is then expressed by a Hamiltonian (total energy)
\bq
H: \mathcal{X} \to \mR,
\eq
defining the vector $e_S$ of as $e_S = \nabla H(x)$, with $\nabla H(x)$ is the column vector of partial derivatives of $H$. Obviously, this implies 
\bq
\frac{d}{dt}H(x(t))= \left(\nabla H(x(t))\right)^T \dot{x}(t)= - e_S^T(t)f_S(t),
\eq 
i.e., the increase of stored energy is equal to the power delivered to the energy-storing elements through the left link in Figure \ref{fig:pHsystems}.

Furthermore, energy dissipation is any relation $\mathcal{R}$ between the flow and effort variables $f_R,e_R$ of the energy-dissipating port having the property that
\bq
\label{R}
e_R^Tf_R \leq 0, \quad (f_R,e_R) \in \mathcal{R}
\eq
Consider now a Dirac structure $\D$ on the manifold $\X \times \mathcal{F}_R \times \mathcal{F}_P$, which is independent\footnote{This can be formalized as a {\it symmetry} of the Dirac structure: the Dirac structure $\D$ is {\it invariant} with respect to arbitrary transformations on $\mathcal{F}_R \times \mathcal{F}_P$; see \cite{vdssym}.} of the point in $\mathcal{F}_R \times \mathcal{F}_P$; i.e., for every $x \in \mathcal{X}$
\bq
\mathcal{D}(x) \subset T_x \mathcal{X} \times \mathcal{F}_R \times \mathcal{F}_P \times T_x^*\mathcal{X} \times \mathcal{E}_R \times \mathcal{E}_P
\eq
is a Dirac structure as in Definition \ref{def:schaft_dn2.1}. 

Then the triple $(\D,H,\R)$, with the energy storage $H: \mathcal{X} \to \mR$, and energy dissipation $\mathcal{R} \subset \mathcal{F}_R \times \mathcal{E}_R$ defines a {\it port-Hamiltonian system} (sometimes abbreviated as {\it pH system}), geometrically defined as the implicit dynamics
\bq 
\label{id}
\begin{array}{c}
(-\dot{x}(t),f_R(t),f_P(t), \nabla H(x(t)), e_R(t), e_P(t)) \in \mathcal{D}(x(t))\\[2mm]
(f_R(t),e_R(t)) \in \mathcal{R}, \quad t \in \mR
\end{array}
\eq
in the state variables $x$, with external port-variables $f_P,e_P$. 

A specific class of port-Hamiltonian systems is obtained by considering Dirac structures which are the graph of a skew-symmetric bundle map
\bq
\begin{bmatrix} - J(x) & -G_R(x) & -G(x) \\ G_R^T(x) & 0 & 0 \\G(x) & 0 & 0 \end{bmatrix}, \quad J(x)=-J^T(x), \quad x \in \X ,
\eq
from $e_S,e_R,e_P$ to $f_S,f_R,f_P$, and a linear energy dissipation relation $e_R=-\bar{R}e_R$ for some matrix $\bar{R}=\bar{R}^T\geq 0$. This yields {\it input-state-output} port-Hamiltonian systems
\bq
\label{iso}
\begin{array}{rcl}
\dot{x} & = & \left[ J(x) - R(x) \right] \nabla H(x) + G(x)u \\[2mm]
y & = & G^T(x) \nabla H(x),
\end{array}
\eq
where $R(x)=G_R(x)\bar{R}\,G_T(x)$, and $u=e_P$ is the {\it input} and $y=f_P$ the {\it output} vector. This is often taken as the starting point for control purposes \cite{passivitybook}. 

On the other hand, for a general Dirac structure {\it algebraic constraints} in the state variables $x$ may easily appear; leading to port-Hamiltonian systems which are {\it not} of the form \eqref{iso}. In fact, if the projection $\rho^*(x)(\D(x))$ of $\D(x)$ to $T^*_x\X$ under the canonical projection $\rho^*(x): T_x \X \times T_x^*\X \to T_x^*\X$ is a strict subspace of $T^*_x \X$, then necessarily $x$ should be such that $\nabla H(x) \in \rho^*(x)(\D(x))$; leading to algebraic constraints \cite{phDAE}, \cite{dalsmo}, \cite{passivitybook}. In the sequel these algebraic constraints, stemming directly from port-based modeling, will be referred to as {\it Dirac algebraic constraints}.

\section{Port-Hamiltonian systems with implicit energy storage and Lagrange algebraic constraints}
An interesting extension of the standard nonlinear port-Hamiltonian systems as discussed in the previous section is obtained as follows. 

For any Hamiltonian $H: \X \to \mR$ the submanifold 
\bq
\mathrm{graph\,}  \nabla H := \{(x,\nabla H(x) ) \mid x \in \X \}
\eq
is a {\it Lagrangian submanifold} \cite{arnold} of the cotangent bundle $T^*\X$. Thus, instead of considering energy storage defined by a Hamiltonian $H$ we may as well consider a general {\it implicit} energy storage defined by a general Lagrangian submanifold $\L$. In fact \cite{arnold} a general Lagrangian submanifold $\L \subset T^*\X$ will be of the form $\mathrm{graph\,} \nabla H$ for a certain $H$ if and only if the {\it projection} of $\L \subset T^*\X$ to $\X$ under the canonical projection $\pi: T^*\X \to \X$ is equal to the whole of $\X$. On the other hand, if and only if the projection $\pi(\L)$ of $\L \subset T^*\X$ to $\X$ is {\it not} equal to the whole of $\X$, then a new class of algebraic constraints arises, namely $x \in \pi(\L)$. These algebraic constraints will be called {\it Lagrange algebraic constraints}; extending the terminology in the linear case in \cite{DAE}. 

The resulting triple $(\D,\L,\R)$ will be called a {\it generalized} port-Hamiltonian system, defining the dynamics (generalizing \eqref{id})
\bq 
\begin{array}{c}
(-\dot{x}(t),f_R(t),f_P(t), e_S(t), e_R(t), e_P(t)) \in \mathcal{D}(x(t))\\[2mm]
(f_R(t),e_R(t)) \in \mathcal{R}, \quad (x(t),e_S(t)) \in \L, \quad t \in \mR
\end{array}
\eq
in the state variables $x$, with external port-variables $f_P,e_P$.

The basic idea of this definition (without the inclusion of an energy dissipation relation and external port) can be already found in \cite{barbero}. The definition of the generalized port-Hamiltonian system $(\D,\L,\R)$ extends the definition in the {\it linear} case as recently given in \cite{DAE}; partly motivated by \cite{beattie}.

\subsection{Properties of the Legendre transform}
Before going on with a discussion of the properties of generalized port-Hamilto-nian systems $(\D,\L,\R)$ and their Lagrange algebraic constraints, let us recall the basic properties of the Legendre transform.

Consider a smooth function $P: \X \to \mathbb{R}$, with column vector of partial derivatives denoted by $\nabla P(x)$.
 The {\it Legendre transform} of $P$ is defined in local coordinates $x$ for $\X$ as the expression
 \bq
 \label{leg}
 P^*(e)= e^Tx - P(x), \quad e = \nabla P(x),
 \eq
 where $e$ are corresponding coordinates for the cotangent space. In the expression \eqref{leg} it is meant that $x$ is expressed as a function of $e$ through the equation $e = \nabla P(x)$; thus obtaining a function $P^*$ of $e$ only. This requires that the map $x \mapsto  \nabla P(x)$ is injective\footnote{However, more generally, i.e., without this assumption, we can still define $P^*$ as the restriction of the function $e^Tx - P(x)$ defined on $T^*\X$, but {\it restricted} to the submanifold $e = \nabla P(x)$. On this submanifold obviously the partial derivatives of $e^Tx - P(x)$ with respect to $x$ are zero, and thus the function is determined as a function of $e$ only.}.
 
 A well-known property of the Legendre transform is the fact that the Legendre transform of $P^*$ is equal to $P$; i.e., $P^{**}=P$. Furthermore, the inverse of the map $ x \mapsto e = \nabla P(x)$ is given as $ e \mapsto x = \nabla P^*(e)$, that is
 \bq
 \nabla P^*(\nabla P(x))=x, \quad \nabla P(\nabla P^*(e))=e
 \eq
Given the Legendre transform $P^*$ of $P$ one may also define the new function
\bq
\widetilde{P}(x) := P^*(\nabla P(x))
\eq
In case of a {\it quadratic} function $P(x)=x^TMx$ for some invertible symmetric matrix $M$ it is obvious to check that $\widetilde{P}=P$; however for a general $P$ this need not be the case. 

Interestingly, using $P^*(e)= e^Tx - P(x)$, $x = \nabla P^*(e)$, and the identity $\nabla P(\nabla P^*(e))=e$, the function $\widetilde{P}$ can be also expressed as
\bq
\widetilde{P}(x) = (\nabla P(x))^T \nabla P^*(\nabla P(x)) - P(\nabla P^*(\nabla P(x))) = x^T \nabla P(x) - P(x)
\eq
Furthermore, we note the following remarkable property
\bq
\nabla \widetilde{P}(x) = \nabla^2 P(x) \nabla P^*(\nabla P(x)) = \nabla^2 P(x)x,
\eq
with $\nabla^2 P(x)$ denoting the Hessian matrix of $P$.

Finally, all of this theory can be extended to {\it partial} Legendre transformations. Consider a partitioning $I \cup J= \{1,\cdots,n\}$ , and the corresponding splitting $x=(x_I,x_J)$, $e=(e_I,e_J)$. The partial Legendre of $P(x_I,x_J)$ with respect to $x_J$ is defined as
\bq
P^*(x_I,e_J) = e_J^Tx_J - P(x_I,x_J), \quad e_J= \frac{P}{\partial x_J}(x_I,x_J),
\eq
where $x_J$ is solved from $e_J= \frac{\partial P}{\partial x_J}(x_I,x_J)$.

\subsection{Explicit representation of implicit storage relations}
In this subsection we will show how generalized port-Hamiltonian systems $(\D,\L,\R)$ with implicit energy storage relations $\L \subset T^*\X$ can be explicitly represented as follows. This extends the observations made in the linear case \cite{beattie}, \cite{DAE} in a non-trivial way.

The starting point is the fact that any Lagrangian submanifold $\L \subset T^*\X$, with $\dim \X =n$, can be locally written as \cite{arnold}
\bq
\label{L}
\L = \{(x,e_S)=(x_I,x_J,e_I,e_J) \in T^*\X \mid e_I = \frac{\partial V}{\partial x_I}, x_J=-\frac{\partial V}{\partial e_J} \} ,
\eq
for some splitting $\{1, \cdots,n\}=I \cup J$ of the index set, and a function $V(x_I,e_J)$, called the {\it generating function} of the Lagrangian submanifold $\L$. In particular, $x_I,e_J$ define local coordinates for $\L$. Now define the Hamiltonian $\widetilde{H} (x_I,e_J)$ as
\bq
\widetilde{H} (x_I,e_J):= V(x_I,e_J) - e^T_J \frac{\partial V}{\partial e_J}(x_I,e_J)
\eq
By Equation \ref{L} the coordinate expressions of $f_S=-\dot{x},e_S$ (in terms of $x_I,e_J$) are given as
\bq
-f_S= \begin{bmatrix} I & 0 \\ -\frac{\partial^2 V}{\partial e_J \partial x_I} & -\frac{\partial^2 V}{\partial e^2_J} \end{bmatrix}
\begin{bmatrix} \dot{x}_I \\ \dot{e}_J \end{bmatrix}, \quad e_S = \begin{bmatrix} \frac{\partial V}{\partial x_I} \\ e_J \end{bmatrix} 
\eq
This yields 
\bq
\frac{d}{dt} \widetilde{H}(x_I,e_J)= --e_S^T(t)f_S(t)
\eq
Consider furthermore any modulated Dirac structure $\D(x) \subset T_x\X \times T_x^*\X \times \F_R \times \E_R \times \F_P \times \E_P$. Since by the power-conservation property of Dirac structures $-e_S^Tf_S = e_R^Tf_R + e_P^Tf_P$ it thus follows that
\bq
\frac{d}{dt} \widetilde{H}(x_I,e_J)= e_R^T(t)f_R(t) + e_P^T(t)f_P(t) \leq e_P^T(t)f_P(t)
\eq
Hence $\widetilde{H} (x_I,e_J)$ serves as the expression of an {\it energy function} (however, {\it not} in the original state variables $x$, but instead in the mixed set of coordinates $x_I,e_J$).

Note that if the mapping $x_J=-\frac{\partial V}{\partial e_J}$ from $e_J$ to $x_J$ is {\it invertible}, and hence the Lagrangian submanifold can be also parametrized by $x=(x_I,x_J)$, and thus $\L=\{(x,\nabla H(x)) \mid x \in \X \}$ for a certain $H$, then actually $\widetilde{H}(x_I,x_J)$ is minus the partial Legendre transform of $V(x_I,e_J)$ with respect to $e_J$, and equals $H$.

As another special case, let us take $I$ to be empty and thus $e_J=e_S$. Let $V^*(x)$ be the Legendre transform of $V(e_S)$. Then it follows that
\bq
V(e_S) - e_S^T \nabla V(e_S) = - \widetilde{V}(e_S)
\eq

\section{Transformation of Dirac algebraic constraints into Lagrange algebraic constraints, and back}
In the previous sections it was discussed how in standard port-Hamiltonian systems $(\D,H,\R)$ {\it Dirac algebraic constraints} may arise (whenever the projection of the Dirac structure onto the cotangent space of the state space is not full), while generalized port-Hamiltonian systems $(\D,\L,\R)$ may have (additional) {\it Lagrange algebraic constraints} (due to the projection of $\L$ on the state space manifold $\X$ not being surjective).

In this section we will show how one can {\it convert} Dirac algebraic constraints (as favored by port-based modeling) into Lagrange algebraic constraints (sometimes having advantages from a numerical simulation point of view); by adding {\it extra state variables}. This extends the constructions detailed in \cite{DAE} from the linear to the nonlinear case. 

\subsection{From Dirac algebraic constraints to Lagrange algebraic constraints}
The first observation \cite{dalsmo} to be made is that a general Dirac structure $\D$ can be written as the graph of a skew-symmetric map on an extended state space as follows. In fact, suppose $\pi^*(\D(x)) \subset T_x\X^*$ is $(n-k)$-dimensional. Define $\Lambda:= \mR^k$. Then there exists a full-rank $n \times k$ matrix $B(x)$ and a skew-symmetric $n \times n$ matrix $J(x)$ such that
\bq
\label{Diracconstraint}
\D(x) \! = \! \{(f,e) \! \in T_x\X \times T_x^*\X \mid \exists \lambda^* \in \Lambda^* \mbox{ s.t. }-f = J(x)e + B(x)\lambda^*, 0  =  B^T(x)e \}
\eq
Conversely, any such equations for a skew-symmetric map $J(x): T_x^*\X \to T_x\X$ define a Dirac structure. Now, let the energy-storage relation of the port-Hamiltonian system be given in a standard way; i.e., by a Hamiltonian $H: \X \to \mR$. Then with respect to the {\it extended} state space $\X_e:= \X \times \Lambda$ we may define the {\it implicit} energy storage relation given by the Lagrangian submanifold (of the same type as in \eqref{L})
\bq
\L_e := \{ (x,\lambda, e, \lambda^*) \in T^* \X_e \mid e= \nabla H(x), \lambda = 0 \},
\eq
corresponding to the Lagrange algebraic constraint $0=\lambda (\,= B^T(x)\nabla H(x))$. Hence the Dirac algebraic constraint $0=B^T(x) \nabla H(x)$ has been transformed into the Lagrange algebraic constraint $\lambda=0$ on the extended state space $\X_e$. The generating function of $\L_e$ is $H$, which is independent of $\lambda^*$, and therefore $\widetilde{H}(x, \lambda^*):= H(x) - \lambda^{*T}\frac{\partial H}{\partial \lambda^*}(x)=H(x)$. 

%Note that the above transformation of Dirac algebraic constraints into Lagrange algebraic constraints by extension of the state space is opposite to the situation considered in Example \ref{singular}, where Dirac algebraic constraints were transformed into Lagrange algebraic constraints by {\it reduction} (leaving out the $u$ variables), or said differently, where Lagrange algebraic constraints were transformed into Dirac algebraic constraints by extension of the state space. See \cite{DAE} for further information.

\subsection{From Lagrange algebraic constraints to Dirac algebraic constraints}
The conversion of Lagrange algebraic constraints into Dirac algebraic constraints is also based on an extension of the state space. It is based on the fact, see e.g. \cite{barbero} and the references quoted in there, that any Lagrangian submanifold $\L \subset T^*\X$ can be locally represented by a parametrized family of generating functions, called a {\it Morse family}. 

To be precise, given a Lagrangian submanifold $\L \subset T^*\X$, a point $P \in \L$ and projection $\pi(P) \in \X$, there exists a neighborhood $V$ of $\pi(P)$, a natural number $k$, a neighborhood $W$ of $0 \in \mR^k$, together with a smooth function $F: V \times W \to \mR$, such that the rank of $\frac{\partial F}{\partial \lambda}$, with $\lambda \in \mR^k$, is equal to $k$ on $\left(\frac{\partial F}{\partial \lambda}\right)^{-1}(0)$, and 
\bq
\{ \left(x, \frac{\partial F}{\partial x}(x, \lambda) \right) \mid \frac{\partial F}{\partial \lambda}(x, \lambda)=0 \}
\eq
is a neighborhood of the point $P$ in $\L$. The function $F(x,\lambda)$, seen as a function of $x$, parametrized by $\lambda$, is called a Morse family for the Lagrangian submanifold $\L$.

Furthermore, given any (modulated) Dirac structure $\mathcal{D}(x) \subset T_x \mathcal{X} \times \mathcal{F}_R \times \mathcal{F}_P \times T_x^*\mathcal{X} \times \mathcal{E}_R$ as before, one may take the direct product with the (trivial) Dirac structure $\{(f_{\lambda},e_{\lambda}) \mid e_{\lambda}=0 \}$, so as to obtain an extended Dirac structure $\D_e$. This defines an {\it extended} port-Hamiltonian system $(\D_e,F, \R)$ with {\it explicit} energy function $F(x,\lambda)$, and thus without Lagrange algebraic constraints.

\begin{example}[Optimal control \cite{barbero}]\label{singular} Consider the optimal control problem of minimizing a cost functional $\int  L(q,u) dt$ for the control system $\dot{q}=f(q,u)$, with $q \in \mR^n, u \in \mR^m$. Define the optimal control Hamiltonian 
\bq
K(q,p,u)=  p^T f(q,u) + L(q,u) ,
\eq
with $p \in \mR^n$ the co-state vector. Application of Pontryagin's Maximum principle leads to the consideration of the standard port-Hamiltonian system (without inputs and outputs) on the space with coordinates $(q,p,u)$, given as
\bq
\label{o1}
\begin{bmatrix}
\dot{q} \\[2mm] \dot{p} \\[2mm] 0 
\end{bmatrix} = 
\begin{bmatrix}
0 & I_n & 0 \\[2mm]
-I_n & 0 & 0 \\[2mm]
0 & 0 & I_m 
\end{bmatrix} 
\begin{bmatrix}
\frac{\partial H}{\partial q}(q,p,u) \\[2mm]
\frac{\partial H}{\partial p}(q,p,u) \\[2mm]
\frac{\partial H}{\partial u}(q,p,u) 
\end{bmatrix}
\eq
The underlying Dirac structure is given as
\bq
\D =\{ \left( \begin{bmatrix} f_q \\[2mm] f_p \\[2mm] f_u \end{bmatrix} , \begin{bmatrix} e_q \\[2mm] e_p \\[2mm] e_u \end{bmatrix}\right) \mid f_q = - e_p, \, f_p = e_q, \, e_u=0\} ,
\eq
i.e., the direct product of the Dirac structure on the $(q,p)$ space given by the graph of the canonical skew-symmetric map $\begin{bmatrix} 0 & -I \\ I & 0 \end{bmatrix}$ from $\begin{bmatrix} f_q \\ f_p \end{bmatrix}$ to $\begin{bmatrix} e_q \\ e_p \end{bmatrix}$, with the trivial Dirac structure $\{(f_u,e_u) \mid e_u=0 \}$. The resulting Dirac algebraic constraint is $\frac{\partial H}{\partial u}(q,p,u) =0$.

System \eqref{o1} can be equivalently rewritten as a port-Hamiltonian system system only involving the $(q,p)$ variables, with {\it implicit} energy storage relation given by the Lagrangian submanifold
\bq
\label{o2}
\L = \{\left( 
\begin{bmatrix}
q \\[2mm] p 
\end{bmatrix}, 
\begin{bmatrix}
e_q \\[2mm] e_p 
\end{bmatrix} \right) 
\mid 
\exists u \mbox{ s.t. } 
\begin{bmatrix}
e_q \\[2mm] e_p 
\end{bmatrix}  = 
\begin{bmatrix}
\frac{\partial H}{\partial q}(q,p,u) \\[2mm]
\frac{\partial H}{\partial p}(q,p,u) 
\end{bmatrix} , \; \frac{\partial H}{\partial u}(q,p,u)=0
\}
\eq
Thus the function $H(q,p,u)$ defines a Morse family (function of $(q,p)$ parametri-zed by $u$) for this Lagrangian submanifold, and the conversion of \eqref{o2} into \eqref{o1} is an example of conversion of a Lagrange algebraic constraint into a Dirac algebraic constraint. See for the linear case \cite{DAE}.
\end{example}

\section{Conclusions}
We have laid down a framework for studying Dirac and Lagrange algebraic constraint equations as arising in (generalized) port-Hamiltonian systems, extending the linear results of \cite{DAE}, \cite{beattie}. In particular, a definition is provided of a nonlinear generalized port-Hamiltonian system, extending the one in \cite{barbero} by including energy dissipation and external ports. Furthermore, we have shown how implicit energy storage relations locally can be explicitly represented by a Hamiltonian depending on part of the state variables and a complementary part of the co-state variables. Also, by extension of the state space (inclusion of Lagrange multipliers) we have shown how Dirac algebraic constraints can be converted into Lagrange algebraic constraints, and conversely.

This work should be seen as a starting point for further study on the numerical properties of the resulting structured classes of nonlinear DAE systems; including their index and regularization \cite{kunkel}. Also it motivates the development of control theory for classes of physical nonlinear DAE systems, as well as extensions to the distributed-parameter case.

%Text with citations \cite{RefB} and \cite{RefJ}.
%\subsection{Subsection title}
%\label{sec:2}
%as required. Don't forget to give each section
%and subsection a unique label (see Sect.~\ref{sec:1}).
%\paragraph{Paragraph headings} Use paragraph headings as needed.
%\begin{equation}
%\end{equation}
%
%% For one-column wide figures use
%\begin{figure}
%% Use the relevant command to insert your figure file.
%% For example, with the graphicx package use
%  \includegraphics{example.eps}
%\end{figure}
%%
%% For two-column wide figures use
%\begin{figure*}
%% Use the relevant command to insert your figure file.
%% For example, with the graphicx package use
%  \includegraphics[width=0.75\textwidth]{example.eps}
%\end{figure*}
%\begin{table}
%\begin{tabular}{lll}
%\hline\noalign{\smallskip}
%\end{tabular}
%\end{table}

\begin{acknowledgements}
Dedicated to Volker Mehrmann at the occasion of his $65$-th birthday, for many stimulating and inspiring discussions on the modeling, simulation and control of complex physical systems.
\end{acknowledgements}

\section*{Conflict of interest}
The authors declare that they have no conflict of interest.



\end{document}